\documentclass{article}
\usepackage{amssymb}
\usepackage{amsmath}

\usepackage{amssymb}
\usepackage{amsfonts}
\usepackage{amscd}

\newtheorem{theorem}{Theorem}

\newtheorem{corollary}[theorem]{Corollary}

\newtheorem{definition}[theorem]{Definition}

\newtheorem{lemma}[theorem]{Lemma}

\newtheorem{proposition}[theorem]{Proposition}
\newtheorem{remark}[theorem]{Remark}

\newcommand\C{\mathbb{C}}
\newcommand\K{\mathbb{K}}
\newcommand\p{(\mathcal{P},\cdot )}
\newcommand\pp{\mathcal{P}}
\begin{document}

\title{Non-associative algebras associated to Poisson algebras}

\author{Michel GOZE , Elisabeth REMM}

\date{}
\maketitle

{\abstract Poisson algebra is usually defined to be a commutative algebra
together with a Lie bracket, and these operations are required to satisfy
the Leibniz rule. We describe Poisson structures in terms of a single
bilinear operation. This enables us to explore Poisson algebras in the realm of non-associative
algebras. We
study their algebraic and cohomological properties, their deformations  
as non-associative algebras,  and
give a classification in low dimensions.}

\section{Poisson algebras presented as non-associative algebras}

Let $\K$ be a commutative field of characteristic different from $2$ and $3$.

\subsection{Non-associative algebra associated to a Poisson algebra}

A Poisson algebra over $\mathbb{K}$ is a $\mathbb{K}$-vector space $\mathcal{P}$
equipped with two bilinear operations:

\noindent 1) A Lie bracket, referred to as the Poisson bracket, usually
denoted by $\{ \ , \ \}$. 

\noindent 2) An associative commutative multipliction which we denote it by 
$\bullet$.

\smallskip

\noindent These two operations are required to satisfy  Leibniz condition: 
\begin{eqnarray}
 \label{leib}
\{X   \bullet Y,Z \}=X\bullet \{Y,Z \}+\{ X,Z\}\bullet Y,
 \end{eqnarray}
for all $X,Y,Z$ in $\mathfrak{g}.$
 This condition means that, with respect to each of the two variables, the Poisson bracket 
behaves as a derivation relative to the multiplication. We denote the Poisson algebra
 by $(\mathcal{P}, \{,\},\bullet )$.

\medskip

Let $\cdot :(X,Y)\rightarrow X\cdot Y$ be a bilinear map on the $\K$-vector space $\mathcal{P}$. The associator
$A$ of $\cdot $ is the trilinear map on $\mathcal{P}$ given by
$$A(X,Y,Z)=(X\cdot Y)\cdot Z-X\cdot (Y\cdot Z).$$

Throughout the paper we do not assume our algebras to be necessarily associative.

\begin{proposition}
Let $(\mathcal{P},\cdot )$ be a  $\mathbb{K}$-algebra. Let us consider the multiplications
$\{, \}$ and $\bullet$ on $\mathcal{P}$
 given by
\begin{eqnarray}
\{X,Y \} &=&\frac{1}{2}(X\cdot Y-Y\cdot X), \label{poisson bracket}
\\
X \bullet Y & =& \frac{1}{2}(X\cdot Y+Y\cdot X).
\label{poisson bullet}
\end{eqnarray}
\noindent Then $(\mathcal{P},\{,\},\bullet )$ is  a
Poisson algebra  if and only if the multiplication $X\cdot Y$
satisfies:
\begin{eqnarray}
\label{associator} 3A(X,Y,Z)=(X\cdot Z)\cdot Y+(Y\cdot
Z)\cdot X-(Y\cdot X)\cdot Z-(Z\cdot X)\cdot Y.
\end{eqnarray}
\end{proposition}

\noindent {\it Proof:} See \cite{M.R}.

\smallskip

\begin{definition}
A non-associative $\K$-algebra $( \mathcal{P},\cdot )$ whose associator satisfies Equation
\rm{(\ref{associator})} is
called an admissible Poisson algebra.
\end{definition}
Let  $\p$ and $(\mathcal{P},\star )$ be admissible Poisson algebras 
defining the same Poisson algebra $(\mathcal{P},\{,\},\bullet )$. Then 
$$
\begin{array}{l}
X\cdot Y-Y\cdot X=X\star Y-Y\star X=2\{X,Y\}\\
X\cdot Y+Y\cdot X=X\star Y+Y\star X=2X\bullet Y
\end{array}
$$
and $X\cdot Y=X\star Y$ because the characteristic of $\K$ is not $2$.
\begin{proposition}
Every Poisson algebra $(\mathcal{P},\{,\},\bullet )$ is associated to precisely one admissible
Poisson algebra $\p$.
\end{proposition}

\noindent Given a Poisson algebra $(\mathcal{P},\{,\},\bullet )$, we shall say that 
$(\mathcal{P}, \cdot)$ where 
$$X \cdot
Y = \{ X,Y\}+X \bullet Y$$ is the admissible Poisson algebra associated 
$(\mathcal{P},\{,\},\bullet )$.  The corresponding Lie
algebra $(\mathcal{P}, \{,\})$ will be denoted by
$\frak{g}_{\mathcal{P}}$ and the associative commutative
algebra $(\mathcal{P}, \bullet )$  by
$\mathcal{A}_{\mathcal{P}}$.

\medskip

\noindent{\bf Notation. }We will denote
 (when no confusion is possible) the Poisson product
 by $XY$ instead of $X\cdot Y$.

\begin{proposition}
An admissible Poisson algebra $\p$ is flexible, that is, the associator satisfies
$$A(X,Y,X)=0$$
for every $X,Y \in \pp.$
\end{proposition}
{\it Proof. } From (\ref{associator}) we have
$$3A(X,Y,X)=X^2Y+(YX)X-(YX)X-X^2Y=0$$
where $X^2=XX.$ Then $\p$ is flexible. $\Box$

\medskip

\medskip

\noindent We deduce easily that the associator of the multiplication $\cdot $ satisfies
\begin{equation}
\label{flex} A(X,Y,Z)+A(Z,Y,X)=0  \quad
{\rm(flexibility)}
\end{equation}
\begin{equation}
\label{equa}
A(X,Y,Z)+A(Y,Z,X)-A(Y,X,Z)=0.
\end{equation}
This last relation is obtained by writing identity
(\ref{associator}) for the triples $(X,Y,Z)$, $(Y,Z,X)$ and
$(Y,X,Z)$.

\begin{remark} {\rm The  system $\{ (\ref{flex}), (\ref{equa})\}$  
is equivalent to the equation
\begin{eqnarray}
\label{sigma3}
2A(X,Y,Z)+\frac{1}{2}A(Y,X,Z)+ A(Z,Y,X)+A(Y,Z,X) +\frac{3}{2}A(Z,X,Y)=0.
\end{eqnarray}
In fact $(\ref{flex})+(\ref{equa})$ implies $(\ref{sigma3})$. Conversely if $(\ref{sigma3})$ is satisfied,
then $(\ref{sigma3})$ applied to the triple $(X,Y,X)$ gives
$$2A(X,Y,X)+A(Y,X,X)+A(X,X,Y)=0$$
and to the triple $(X,X,Y)$
$$5A(X,X,Y)+5A(Y,X,X)+2A(X,Y,X)=0.$$
We deduce (\ref{flex}) and (\ref{equa}).
It is worth noting that a non-associative algebra satisfying (\ref{sigma3}) is not always an admissible Poisson algebra.}
\end{remark}

\begin{proposition}
An admissible Poisson algebra $\p$ is a power associative algebra.
\end{proposition}
{\it Proof.} Recall that a non-associative algebra is  power
associative if every element generates an associative subalgebra.
Let $X$ be in $\p$. We define the power of $X$ by $X^1=X, \,
X^{i+1}=X\cdot X^i$. We will prove that
$X^{i+n}X^{j-n}=X^{i-p}X^{p}=X^{i+j}$ for all $i,j\geq 1$ and $1
\leq p \leq i, \ 1 \leq n \leq j$.
Since $\p$ is flexible, we have  $A(X,X^{j},X)=0$ for any $1\leq j$. We have  
$X^{j}X=XX^j$ for $j=1$. Suppose that this equation is true for $j$, then 
$A(X,X^j,X)=0$ and $X^{j+1}X=X(X^jX)=X(XX^j)=XX^{j+1}$. So for any $j \geq 1, \
X^jX=XX^j$. Now we shall use induction over $i$ to prove
that, for any $j\geq 1$, $X^iX^j=X^jX^i$. This identity is
trivial for $i=1$. Suppose that it is satisfied for  $i\geq 1$. 
Then  relation (\ref{associator}) gives
$$3A(X,X^{i},X^j)-(XX^j)X^i-(X^iX^j)X+(X^iX)X^j+(X^jX)X^i= 0$$
and as $X^iX^j=X^jX^i$, we obtain
$$4X^{i+1}X^j= 3X(X^{i}X^j)+(X^{i}X^j)X.$$ 
Similarly, (\ref{associator}) applied to the triple
$(X,X^j,X^i)$ gives
$$4X^{j+1}X^i= 3X(X^{j}X^i)+(X^{j}X^i)X.$$
By assumption $X^iX^j=X^jX^i$, we obtain $X^iX^{j+1}=X^{j}X^{i+1}$. 
By (\ref{associator}),
this implies
 $A(X^i,X,X^j)=0$. Thus,
$$X^{i+1}X^j=X^iX^{j+1}=X^{j}X^{i+1}$$
and $X^iX^j=X^jX^i$ for all $i, j$. Finally, we prove that for
any $i$ the relation $X^{i-p} X^{p}=X^i$ is satisfied for any $1 \leq
p < i$. It is evident for $i=1$. 
Suppose that these relations are satisfied for a fixed $i$.
Then 
$$3A(X^{i-p},X,X^p)= X^{p+1}X^{i-p} -X^{i-p+1}X^p $$
implies $X^{i-p+1}X^{p} =X^{i-p}X^{p+1} $ and
$$3A(X^{i-p},X^p,X)= X^{p+1}X^{i-p} -X^{i+1} $$
implies  $X^{i+1}=X^{p+1}X^{i-p}$. Thus $X^{i+1-p} X^{p}=X^{i+1}$ 
and the algebra $\p$ is power associative. $\Box$

\medskip

\begin{remark}{\it Poisson algebras as $\mathbb{K}\,[\Sigma_3]$-associative algebras.}

\rm{ \noindent In \cite{G.R2}, large classes of non-associative algebras were
studied. In this section we show that admissible Poisson algebras belong to
this category of algebras.

\noindent Let $\Sigma_3$ be the order three symmetric group and $\mathbb{K}\,[\Sigma_3]$ its $\mathbb{K}$-group
algebra. A (non-associative) $\mathbb{K}$-algebra $(\mathcal{A},\mu)$ is called a
$\mathbb{K}\,[\Sigma_3]$-associative algebra if there exists $v
\in \mathbb{K}\,[\Sigma_3]$, $v \neq 0$, such that
$$A_{\mu} \circ  \Phi_v=0,$$
where $A_{\mu}=\mu\circ (\mu \otimes Id)-\mu \circ (Id\otimes \mu )$ is the associator of the algebra
$\mathcal{A}$ and $\Phi_v: \mathcal{A}^{\otimes 3} \rightarrow
\mathcal{A}^{\otimes 3} $ is defined by
$$\Phi_{\sigma}(v_1,v_2,v_3)=(v_{\sigma^{-1}(1)},v_{\sigma^{-1}(2)},v_{\sigma^{-1}(3)})$$
for all $\sigma \in \Sigma_3$.

\noindent Now suppose that $\p$ is an admissible Poisson algebra. From~(\ref{associator}) we see
that the associator of the multiplication satisfies
$$A_{\mu} \circ  \Phi_{v_1}=0$$
for $v_1=Id-\tau_{12}+c_1$, where $\tau_{ij}$ interchanges elements i and j and 
$c_1(1,2,3)=(2,3,1)$.
The flexibility identity~(\ref{flex}) can be written as $A_{\mu}
\circ  \Phi_{v_2}=0$ for $v_2=Id+\tau_{13}$. Recalling  the
classification of \cite{G.R2}, we deduce that any Poisson algebra
is an algebra of type $(IV_1)$ for $\alpha=-\frac{1}{2}$ (we have
$v=2Id+ \frac{1}{2}\tau_{12}+\tau_{13} +c_1+  \frac{3}{2}c_2$ and
$F_v$ is 4-dimensional).
}
\end{remark}

\subsection{Pierce decomposition}

We say that a power associative algebra $\mathcal{P}$ is a nilalgebra if any element $X$ 
is nilpotent, i.e.
$$\forall X\in \mathcal{P}, \exists r\in \mathbb{N}\  \mbox{\rm such that} \
X^r=0.$$

\begin{proposition}
\label{nie}
Any finite dimensional admissible Poisson algebra which is not a nilalgebra contains
 a non-zero idempotent element.
\end{proposition}
This is a consequence of the power associativity of a Poisson algebra.

\medskip

Let $e$ be a non-zero idempotent, i.e. $e^2=e$. Equation~(\ref{poisson bullet})
implies $e \bullet e =e$,
thus $e$ is an idempotent of the associative algebra $\mathcal{A}_{\mathcal{P}}$.
The Leibniz identity implies
$$\{ e,x \}=\{ e\bullet e,x \}=2e\bullet \{ e,x \}.$$
Therefore, $\{ e,x \}$ is either zero or an eigenvector of the
operator
$$L^{\bullet}_e:x \rightarrow e\bullet x$$
in $\mathcal{A}_{\mathcal{P}}$ associated to the eigenvalue
$\frac{1}{2}.$ Since $e$ is an idempotent, the eigenvalues
associated to $L^{\bullet}_e$ are $1$ or $0$. It follows that $\{
e,x \}=0$ which implies that $e \in Z(\frak{g}_\mathcal{P})$ and
$e \cdot x=e \bullet x=x \bullet e=x \cdot  e$.

\begin{proposition}
Let $\p$ be an admissible Poisson algebra such that the center of the associated Lie algebra $\mathfrak{g}_{\mathcal{P}}$
is  zero. Then $\p$ has no idempotent different from zero. If $\mathcal{P}$
is of finite dimension then it is a nilalgebra.
\end{proposition}

Suppose that there exists an idempotent $e\neq 0.$ Since
$\mathcal{P}$ is flexible, the operators $L^{\bullet}_e$ and
$R^{\bullet}_e$ defined by $L^{\bullet}_e(x)=e \bullet x$ and
$R^{\bullet}_e(x)=x \bullet e$ commute and
$L_e^{\bullet}=L_e^{\cdot},R_e^{\bullet}=R_e^{\cdot}$. Then
$\mathcal{P}$ decomposes as
$$\mathcal{P}=\mathcal{P}_{0,0} \oplus \mathcal{P}_{0,1} \oplus \mathcal{P}_{1,0}\oplus \mathcal{P}_{1,1}$$
with $\mathcal{P}_{i,j}=\left\{x_{i,j}\in \mathcal{P} \ {\mbox{\rm
such \ that}} \ ex_{i,j}=ix_{i,j},\, x_{i,j}e=jx_{i,j} \right\},$
$i,j \in \left\{ 0,1 \right\}.$ From Proposition \ref{nie}, $e
\in  Z(\mathfrak{g}_{\mathcal{P}}).$ So $\{e,x\}=0$ for any $x$,
that is,  $ex=xe$ and $
\mathcal{P}_{0,1}=\mathcal{P}_{1,0}=\left\{0\right\}.$

\begin{proposition}
If the admissible Poisson algebra $\p$ has a non-zero idempotent, it admits the Pierce decomposition
$$\mathcal{P}=\mathcal{P}_{0,0} \oplus \mathcal{P}_{1,1},$$
where $\mathcal{P}_{0,0}$ and $ \mathcal{P}_{1,1}$ are admissible Poisson algebras with the induced product.
\end{proposition}

\noindent {\it Proof.} We have to show that $\mathcal{P}_{0,0}$
and $ \mathcal{P}_{1,1}$ are Poisson subalgebras. Let $x,y \in
\mathcal{P}_{0,0},$ then $ex=ey=xe=ye=0.$ From~(\ref{associator}),
we obtain
$$
\left\{
\begin{array}{lll}
 -3e(xy) & = & (xy)e \\
0 & = & (xy)e-(yx)e \\
3(xy)e & = & -(yx)e.
\end{array}
\right.
$$
So $(xy)e=-3e(xy)=(yx)e=-3(xy)e$ and $(xy)e=e(xy)=0$. Then $xy \in \mathcal{P}_{0,0}.$ Similarly
if $x,y \in \mathcal{P}_{1,1}$, then (\ref{associator}) applied to the triple $(e,x,y)$ gives $xy=e(xy)$.
The same equation applied to $(x,e,y)$ and $(x,y,e)$ gives
$$
\left\{
\begin{array}{lll}
 (xy)e+yx-xy-(yx)e & = & 0 \\
3(xy)e -3xy-yx+(yx)e & = & 0.
\end{array}
\right.
$$
Thus, $4(xy)e -4xy=0$ which means that $(xy)e=xy$ and $\mathcal{P}_{1,1}$ is a Poisson subalgebra of $\p$ $\Box$

\begin{remark} \rm{ Poisson algebras are Lie-admissible power-associative algebras. 
In \cite {Ko} Kosier
gave examples of simple Lie-admissible power-associative finite-dimensional algebras called anti-flexible algebras.
These algebras also have the property  $A=A_{00} \oplus A_{11}$ in every Pierce decomposition.}
\end{remark}

\subsection{Pierce decomposition associated to orthogonal idempotents}

Let $e_1$ and $e_2$ be non-zero orthogonal idempotents,  $e_1  e_2=e_2  e_1=0.$ Let
$\mathcal{P}=\mathcal{P}^1_{0,0}\oplus
\mathcal{P}^1_{1,1}=\mathcal{P}^2_{0,0}\oplus \mathcal{P}^2_{1,1}$
be the corresponding Pierce  decompositions. Let us suppose that
$x \in \mathcal{P}^1_{0,0}$. Applying (\ref{associator}) to
the triples associated to the elements $\{e_1,e_2,x\}$, we obtain
the condition
$$(xe_2)e_1=(e_2x)e_1=e_1(e_2x)=e_1(xe_2)=0$$
for the elements $xe_2$ and $e_2x$  in $\mathcal{P}^1_{0,0}$. In other words,
$$L_{e_2}(\mathcal{P}^1_{0,0})\subset \mathcal{P}^1_{0,0}, \ \ R_{e_2}(\mathcal{P}^1_{0,0})
\subset \mathcal{P}^1_{0,0},$$ where $L_{e_2}(x)=e_2x$ and
$R_{e_2}x=xe_2$. So, $e_2$ is an idempotent of the Poisson
algebra $(\mathcal{P}^1_{0,0},.)$. Thus we have
$$\mathcal{P}^1_{0,0}=\mathcal{P}^1_{0,0}\cap \mathcal{P}^2_{0,0}\oplus \mathcal{P}^1_{0,0}\cap \mathcal{P}^2_{1,1}.$$
Using the same reasonings, we can show that if $x \in \mathcal{P}^1_{1,1}$ then, $e_2x=xe_2=0$ and
$$\mathcal{P}^1_{1,1}\subset \mathcal{P}^2_{0,0}.$$
Thus, $\mathcal{P}^1_{0,0}=\mathcal{P}^1_{0,0}\cap
\mathcal{P}^2_{0,0}\oplus \mathcal{P}^2_{1,1}.$ Observe that
$\mathcal{P}^1_{1,1}$ cannot be further decomposed using the
spaces $\mathcal{P}^2_{0,0}$ and $ \mathcal{P}^2_{1,1}$ associated
to $e_2$ as we have 
$$ \mathcal{P}^1_{1,1} =\mathcal{P}^1_{1,1}\cap
\mathcal{P}^2_{0,0}\oplus \mathcal{P}^1_{1,1} \cap
\mathcal{P}^2_{1,1}. $$ But $\mathcal{P}^2_{1,1} \subset
\mathcal{P}^1_{0,0}$ so that $\mathcal{P}^1_{1,1} \cap
\mathcal{P}^2_{1,1}=  \left\{ 0 \right\}$ and $\mathcal{P}^1_{1,1}
\cap \mathcal{P}^2_{0,0}=\mathcal{P}^1_{1,1} $. Then,
$$\mathcal{P}=\mathcal{P}^1_{0,0} \cap \mathcal{P}^2_{0,0}  \oplus
\mathcal{P}^1_{1,1} \oplus \mathcal{P}^2_{1,1}.$$

\begin{proposition}\label{decomp}
If $e_1$ and $e_2$ are non-zero orthogonal idempotents, then $\mathcal{P}$ decomposes into a direct
sum of Poisson subalgebras,
$$\mathcal{P}=\mathcal{P}^1_{0,0} \cap \mathcal{P}^2_{0,0}  \oplus \mathcal{P}^1_{1,1} \oplus \mathcal{P}^2_{1,1}.$$
\end{proposition}
Proposition~\ref{decomp} can be easily generalized to a family of orthogonal
idempotents $\left\{e_1,...,e_k \right\}$. The corresponding
decomposition can then be written as
$$\mathcal{P}=\cap_{i=1}^k \mathcal{P}^i_{0,0}\oplus_{i=1}^k \mathcal{P}^j_{1,1} .$$

\subsection{Radical of a Poisson algebra}

We already know that a Poisson algebra $\p$ is power associative.
Recall that an element $x\in \mathcal{P}$ is nilpotent if there is an integer $r$
such that $x^r=0$. An algebra (two-sided ideal) consisting only of
nilpotent elements is called a nilalgebra (nilideal). If
$\mathcal{P}$ is a finite dimensional Poisson algebra, then there is
a unique maximal nilideal $\mathcal{N}(\mathcal{P})$ called the
nilradical. Let $\mathcal{A}_{\mathcal{P}}$ be the commutative
associative algebra associated to $\p$. Then, the Jacobson radical
$J(\mathcal{A}_{\mathcal{P}})$ of $\mathcal{A}_{\mathcal{P}}$
contains $\mathcal{N}(\mathcal{P}).$ Since
$\mathcal{N}(\mathcal{P})$ is a two-sided ideal of $\p$, it is also a
Lie ideal of $\frak{g}_{\mathcal{P}}.$ One can easily prove:

\begin{proposition}
\label{prop8} The nilradical $\mathcal{N}(\mathcal{P})$ of $\p$
coincides with the maximal Lie ideal of $\frak{g}_\mathcal{P}$
contained in $\mathcal{J}(\mathcal{A}_\mathcal{P}).$
\end{proposition}

\begin{remark} \rm{ \label{nil} In the category of associative algebras, or more
generally, of alternative algebras, any nilalgebra is nilpotent.
This is no longer true in the category of Poisson algebras as the
following example shows.

Let $\p$ be the 3-dimensional algebra defined by
$$
\left\{
\begin{array}{l}
e_i^2=0 \\
e_1e_2=-e_2e_1=e_2 \\
e_1e_3=-e_3e_1=-e_3 \\
e_2e_3=-e_3e_2=e_1.
\end{array}
\right. 
$$
The corresponding  algebra $\mathcal{A}_\mathcal{P}$ is abelian and any element 
of $\mathcal{P}$ is nilpotent.
The Poisson algebra $\mathcal{P}$ is a nilalgebra. But $\mathcal{P}^2=\mathcal{P}$
so $\mathcal{P}$ is not a nilpotent algebra. This algebra is an example of simple nilalgebra.}
\end{remark}

\begin{remark} \rm{ An element $x \in \mathcal{P}$ is {\it properly
nilpotent} if it is  nilpotent and $xy$ and $yx$
are nilpotent for any $y \in \mathcal{P}.$ The Jacobson radical 
of $\mathcal{A}_\mathcal{P}$ coincides with the set of properly
nilpotent elements of $\mathcal{A}_\mathcal{P}.$ Let $x$ be a
properly nilpotent element of $\mathcal{P}$ and suppose that $x
\notin \mathcal{N}(\mathcal{P}).$ We know that $x \in
\mathcal{J}(\mathcal{A}_\mathcal{P}).$ By Proposition
\ref{prop8}, there exists $y \in \mathcal{P}$ such that $\{ x,y\}
\notin \mathcal{N}(\mathcal{P}).$ We have $x \bullet y \in
\mathcal{J}(\mathcal{A}_\mathcal{P})$. This implies that $\{
x,y\}\notin \mathcal{J}(\mathcal{A}_\mathcal{P})$, otherwise $xy
\in \mathcal{J}(\mathcal{A}_\mathcal{P})$ and
$\mathcal{N}(\mathcal{P})$ would not be maximal. But $x \in
\mathcal{J}(\mathcal{A}_\mathcal{P})$, so $xy$ is nilpotent and $xy
\in \mathcal{J}(\mathcal{A}_\mathcal{P})$. This is a
contradiction and the nilradical coincides with the set of
properly nilpotent elements. Zorn's theorem concerning nilalgebra still holds in the
framework of Poisson algebras.}
\end{remark}

\begin{remark} \rm{ We have seen that any finite dimensional Poisson algebra which is not a nilalgebra contains
a non-zero idempotent. An idempotent $e$ is {\it principal} if there is no
idempotent $u$ orthogonal to $e$ (i.e. $ue=eu=0$ with $u^2=u\neq 0$).
If $\p$ is not a nilalgebra, $\mathcal{A}_\mathcal{P}$ is not a nilalgebra
and it has a principal idempotent element. Let $e$ be such an element.
As $e^2=e \bullet e=e$, it is an idempotent element of $\mathcal{P}$.
If one can find $u$ such that $u^2=u \bullet u=u$ with $ue=eu=0$, then $u \bullet e=e \bullet u=0$
which is impossible. Therefore we have:}

\begin{proposition}
Any finite dimensional admissible Poisson algebra which is not a nilalgebra contains a principal idempotent element.
\end{proposition}

\end{remark}

\begin{remark} \rm{ Let us assume that $\mathcal{P}$ is a unitary algebra. If $x$ is an invertible
element of $\mathcal{P}$, there exists $x^{-1} \in \mathcal{P}$ such that $xx^{-1}=x^{-1}x=1$.
In particular $x \bullet x^{-1}=x^{-1} \bullet x=1$ and $x^{-1}$ is the inverse of $x$
in $\mathcal{A}_\mathcal{P}$. Thus the inverse of an invertible element of $\mathcal{P}$
is unique.
Let us note that if $\mathcal{P}$ is unitary, finite dimensional and if the unit is the
only idempotent element, any non-nilpotent element is invertible. In fact,
such an element $x$ generates an associative algebra which admits an idempotent. 
Then $1\in \mathcal{P}$, which turns out to be the only idempotent and can be expressed as
$$1=\sum\alpha_ix^i=x(\sum \alpha_ix^{i-1}).$$
It follows that $\sum \alpha_ix^{i-1}$ is the inverse of $x$.}
\end{remark}

\subsection{Simple Poisson algebras}

An admissible  Poisson algebra $\p$ is {\it simple} if it has not some proper ideal and if
 $\mathcal{P}^2 \neq \{ 0 \}$.
Let $L_x$ and $R_x$ be the left and right translations by  $x \in \mathcal{P}$.
Let $\mathcal{M}(\mathcal{P})$ be the associative subalgebra of $End(\mathcal{P})$ generated
by $L_x,R_x$ for $x \in \mathcal{P}$. In this algebra, we have the following relations
$$
\left\{
\begin{array}{l}
L_x\bullet R_x=R_x\bullet L_x  \\
4L_{x^2}=3(L_x)^2-(R_x)^2+2R_x\bullet L_x \\
4R_{x^2}=3(R_x)^2-(L_x)^2+2R_x\bullet L_x.
\end{array}
\right.
$$

{\it The  algebra $\mathcal{P}$ is simple if and only if $\mathcal{P}$ is a
non-trivial irreducible
$\mathcal{M}(\mathcal{P})$-module.}

One can consider the centralizer
$\tilde{\mathcal{C}}$ of $\mathcal{M}(\mathcal{P})$ in $End(\mathcal{P})$.
If $\mathcal{P}$ is simple and if $\tilde{\mathcal{C}}$ is non-trivial, then
$\tilde{\mathcal{C}}$ is a field which is a central simple Poisson algebra over itself.

\begin{remark} \rm{We saw in Remark \ref{nil} that
there are admissible Poisson algebras which are nilalgebras. In this case
$\mathcal{N}(\mathcal{P})$ is non-zero. We can consider the Albert
radical $\mathcal{R}(\mathcal{P})$ defined as the intersection of
all maximal ideals $\mathcal{M}$ of $\mathcal{P}$ such that
$\mathcal{P}^2 \not\subset \mathcal{M}$. In the algebra defined in Remark~\ref{nil},
$\mathcal{P}^2=\mathcal{P} $. If $\mathcal{M} $ is maximal and
satisfies $\mathcal{M} \subseteq \mathcal{P}^2 $ and $\mathcal{M}
\neq \mathcal{P}^2$, then $\mathcal{M}=\left\{ 0\right\}.$ The
Albert radical is  $\left\{ 0\right\}$ which implies the
semi-simplicity of $\mathcal{P}$.}
\end{remark}

\begin{proposition}
If $\p$ is a simple nilalgebra such that $x^2=0$ for all $x \in \mathcal{P}$
then $\mathcal{A}_\mathcal{P}$ is an associative nilalgebra satisfying
$(\mathcal{A}_\mathcal{P})^2 =0$. 
\end{proposition}

\noindent {\it Proof.} The subalgebra $\mathcal{P}^2=\left\{xy, x,y\in \mathcal{P} \right\}$
  is an ideal of $\mathcal{P},$ so $\mathcal{P}^2=\mathcal{P}$.
By the hypothesis, for every $x \in \mathcal{P}$ we have $x^2=0.$ Then
$$(x+y)^2=x^2+y^2+xy+yx=xy+yx=0$$
for all $x,y \in \mathcal{P}^2$. This implies
$$x \bullet y= \frac{1}{2}(xy+yx)=0$$
thus the associative algebra $\mathcal{A}_\mathcal{P}$ is trivial.

\medskip

We can also consider simple admissible Poisson algebras which are not 

\noindent nilalgebras.
In this case the Albert radical is  $\left\{ 0 \right\}$ and $\mathcal{P}^2 \neq 0.$

\begin{proposition}
Let $\p$ be a finite dimensional simple admissible Poisson algebra which is
not a nilalgebra. Then it  has a unit element.
\end{proposition}

{\it Proof. }In fact $\mathcal{P}$ has a principal idempotent $e$.
Its Pierce decomposition $\mathcal{P}=\mathcal{P}_{0,0}\oplus
\mathcal{P}_{1,1}$ is such that $\mathcal{P}_{0,0} \subset
\mathcal{R}(\mathcal{P})$. Then $\mathcal{P}_{0,0}=\left\{0
\right\}$ and $\mathcal{P}=\mathcal{P}_{1,1}$. Therefore, $e=1$.

\subsection{Classification of simple complex Poisson algebras such
that $\frak{g}_\mathcal{P}$ is  simple }

\begin{lemma}
\label{lemma}\label{lemmasimple}
Let $\p$ be an admissible Poisson algebra.
If $\frak{g}_\mathcal{P}$ is a simple Lie algebra then $\mathcal{P}$ is a simple  algebra.
\end{lemma}
{\it Proof. }If $I\varsubsetneq \mathcal{P}$ is an ideal of $\mathcal{P}$, then $I$ is also an ideal 
of $\frak{g}_\mathcal{P}$
so $I$ must be trivial.

\begin{proposition}
If $\frak{g}_\mathcal{P}$ is a simple complex Lie algebra, then $UV=\{U,V\}$ 
for all $U,V \in \mathcal{P}$, that is,
the associative algebra $\mathcal{A}_\mathcal{P}$ satisfies $\mathcal{A}_\mathcal{P}^2=\{0\}$.
\end{proposition}
{\it Proof. }
Let $\frak{g}_\mathcal{P}$ be a simple complex Lie
algebra of rank $r$. Let $\frak{n}_{-}\oplus \frak{h} \oplus
\frak{n}_{+}$ be its root-decomposition, where $\frak{h}$ is a
Cartan subalgebra. Let $\{Y_j,H_i,X_j\}$ be the corresponding Weyl
basis. Since $\{H_k^2,H_i\}=0$ for all $i=1,...,r$ we deduce that
$$H_k^2 \in \frak{h}, \ \ k=1,...,r.$$
Thus, $H_k^2=\sum_{i=1}^r \alpha _k^iH_i.$ Let us put $\{H_k,X_j\}=\rho
_{k,j}X_j$. We obtain 
$$\{H_k,X_j^2\}=2\rho_{k,j}X_j^2$$
for all $k=1,...,r$. Thus $2\rho_{k,j}$ are also roots of $\frak{g}_\mathcal{P}$, 
but this is impossible
so $X_j^2=0$ for every $j$. Similary we have for all $k=1,...,r$
$$\{H_k,X_j\bullet X_i\}=(\rho_{k,j}+\rho_{k,j})X_j\bullet X_i$$
so $(\rho_{k,j}+\rho_{k,j})$ are roots. This implies
$$X_j\bullet X_i=0.$$
In the same way we have
$$Y_j^0=Y_j\bullet Y_i=0$$
for all $i,j$. It turns out that
$$\{H_k^2,X_j\}=2H_k\bullet \{H_k,X_j\}=2\rho _{k,j}H_k\bullet X_j=\sum_{i=1}^r  \alpha _k^i\rho _{i,j}X_j$$
and
$$\rho _{k,i}H_k\bullet X_j=\frac{1}{2}(\sum_{i=1}^r  \alpha _k^i\rho _{i,j})X_j.$$
For any $j$ there is $k$ such that $\rho _{k,j}\neq 0.$ Thus
$$\{H_k\bullet X_j,X_j\}=0=H_k\bullet \{X_j,x_j\}+\{H_k,X_j\}\bullet X_j=\rho _{k,j}X_j^2$$
and
$$X_j^2=0, \ \forall j.$$
By similar arguments, the identities $Y_j^2=0$ hold. For $i=1,...,r$ we
have 
$$\{X_i^2,Y_i\}=0=2X_i\bullet \{X_i,Y_i\}=-4X_i\bullet H_i.$$
Thus, $\sum \alpha _i^j\rho _{j,i}=0$. As the matrix $(\rho
_{j,i})$ is non-singular, we deduce that $\alpha _i^j=0$, i.e,
$$H_i^2=0, \forall i=1,...,r.$$
The Poisson algebra $\mathcal{P}$ is a nilalgabra. Moreover, $H_i\bullet X_j=H_i\bullet Y_j=0$ 
and we conclude that
$U\bullet V=0$ for all $U,V \in \mathcal{A}_\mathcal{P}$. $\Box$

\section{On the classification of finite dimensional complex Poisson algebras}

Let $\pp$ be a finite dimensional complex Poisson algebra.

\begin{lemma}
\label{lemma2}
If there is a non-zero vector $X \in \frak{g}_{\mathcal{P}}$ such that $ad\, X$ is
diagonalizable
with $0$ as a simple root, then $\mathcal{A}_\mathcal{P}^2=\{0\}$.
\end{lemma}

\noindent {\it Proof.} Let $\left\{ e_1, ..., e_n \right\}$ be a
basis of $\frak{g}_\mathcal{P}$ such that $ad \, e_1$ is diagonal
with respect to this basis. By assumption, $\{e_1,e_i\}= \lambda_i
e_i$ with $\lambda_i \neq 0$ for $i \geq 2$. Since $\{
e_1^2,e_1\}=2e_1 \bullet \{e_1,e_1\}= 0$, it follows that
$e_1^2=ae_1.$ But for any $i \neq 1,$ $\{e_1^2,e_i\}=2e_1 \bullet
\{e_1,e_i\}=2 \lambda_i e_1 \bullet e_i$ and
$\{e_1^2,e_i\}=a\lambda_ie_i$, thus $e_1 \bullet
e_i=\frac{a}{2}e_i$. The associativity of the product $\bullet$
implies that $(e_1 \bullet e_1)\bullet  e_i=ae_1 \bullet
e_i=\frac{a^2}{2}e_i =e_1 \bullet (e_1 \bullet
e_i)=\frac{a^2}{4}e_i$. Therefore $a=0$ and $e_1^2=0=e_1 \bullet
e_i$ for any $i.$ Finally, $0=\{e_1 \bullet e_j,e_i\}=e_1\bullet
\{ e_j,e_i\}+e_j \bullet \{e_1 ,e_i\}=\lambda_i e_j \bullet e_i$,
which implies $ e_i\bullet e_j=0, \ \forall i,j \geq 1.$ $\Box$

\subsection{Classification of 2 dimensional Poisson algebras}

$\bullet$ If $\frak{g}_\mathcal{P}$ is abelian then $\mathcal{A}_\mathcal{P}$ can be any complex
associative commutative algebra and $XY=X\bullet Y$. In this case the classification of Poisson algebras
boils down to the classification of commutative associative complex algebras \cite{Ga}.

\medskip

\noindent $\bullet$ If $\frak{g}_\mathcal{P}$ is not abelian, it is solvable and isomorphic
to the Lie algebra given by $\{e_1,e_2\}=e_2.$ From Lemma~\ref{lemma2} we know that
$\mathcal{A}_\mathcal{P}$ is trivial and $e_ie_j=\{e_i,e_j\}$ for $i,j=1,2.$

\subsection{Classification of 3 dimensional Poisson algebras}

$\bullet$ If $\frak{g}_\mathcal{P}$ is  abelian then $\mathcal{A}_\mathcal{P}$
can be an  arbitrary associative commutative algebra and $XY=X\bullet Y$. 
In this case the classification is given in \cite{Ga}.

 \noindent $\bullet$ If $\frak{g}_\mathcal{P}$ is nilpotent but not abelian it is isomorphic to the
Heisenberg algebra. Let us consider a basis $\left\{ e_i
\right\}_{i=1,2,3}$ of $\frak{g}_\mathcal{P}$ such that
$\{e_1,e_2\}=e_3$. It follows from the Leibniz identities  that
$e_1^2=ae_1+be_3.$ But $\{e_1^2,e_2\}=2e_1 \bullet e_3=ae_3$ and
$\{e_1\bullet e_3,e_2\}=e_3 \bullet e_3=\{ae_3,e_2\}=0.$ The
associativity of $\bullet $ implies that $a=0$. We see that
$$e_1^2=\alpha e_3,  \ e_1 \bullet e_3=e_3^2=0.$$
Similarly,
$$e_2^2=\beta  e_3,  \ e_2 \bullet e_3=0.$$
Finally, $\{e_1\bullet e_2,e_i\}=0$ for $i=1,2,3$ implies
$e_1\bullet e_2=\gamma e_3$. Thus $\mathcal{A}_\mathcal{P}$ is
isomorphic to the algebra:
$$
\left\{
\begin{array}{ll}
e_1^2=\alpha e_3 \\
e_2^2=\beta  e_3 \\
e_1 \bullet e_2=e_2 \bullet e_1=\gamma  e_3. \\
\end{array}
\right. 
$$
We obtain the following Poisson algebra
$$
\left\{
\begin{array}{ll}
e_1^2=\alpha e_3 \\
e_2^2=\beta  e_3 \\
e_1 \cdot e_2=(\gamma  +1)e_3 \\
e_2 \cdot e_1=(\gamma  -1)e_3. \\
\end{array}
\right. 
$$
The base change 
$$
\left\{
\begin{array}{l}
e_1'=ae_1+be_2 \\
e_2'=ce_1+de_2
\end{array}
\right.
$$
gives
$$
\left\{
\begin{array}{l}
(e_1')^2=(a^2\alpha +2ab \gamma  +b^2\beta  )e_3 \\
(e_2')^2=(c^2\alpha +2cd \gamma  +d^2\beta  )e_3.
\end{array}
\right.
$$
If $\gamma  ^2-\alpha \beta  \neq 0,$  the equation $\alpha +2x\gamma  +x^2\beta  =0$ has
two distinct roots and we can assume that $e_1'$ and $e_2'$ are linearly independent such 
that $(e_1')^2=(e_2')^2=0.$
In this case the only possible values of parameters $\alpha $ and $\beta $
are  $\alpha =\beta  =0$. We  obtain the one-parametric family
$$
\mathcal{P}_{3,1}(\gamma  )=
\left\{
\begin{array}{ll}
e_1^2=e_2^2=e_3^2=0 \\
e_1 \cdot e_2=(1+\gamma ) e_3 \\
e_2 \cdot e_1=(-1+\gamma ) e_3 \\
e_1 \cdot e_3=e_3 \cdot e_1=e_3 \cdot e_2=e_2 \cdot e_3=0.
\end{array}
\right. 
$$
If $\gamma  ^2-\alpha \beta  = 0$ and if $\beta  \neq 0$, we can always choose $c$ and $d$ 
such that $e_2^2=0.$ Then we can suppose that $\beta =0$. This implies $\gamma =0$. If
$\alpha =0$ we obtain $
\mathcal{P}_{3,1}(0).$
If $\alpha \neq 0$, we can assume $\alpha =1$ which gives the algebra:
$$
\mathcal{P}_{3,2}=
\left\{
\begin{array}{ll}
e_1^2= e_3\\
e_2^2=e_3^2=0 \\
e_1 \cdot e_2= e_3 \\
e_2 \cdot e_1= - e_3 \\
e_1 \cdot e_3=e_3 \cdot e_1=e_3 \cdot e_2=e_2 \cdot e_3=0.
\end{array}
\right. 
$$
$\bullet$ Suppose that $\frak{g}_\mathcal{P}$ is solvable but not nilpotent. Then
the following three cases may happen.

i) The multiplication is defined by $\{e_1,e_2\}=e_2$. Then $\p$ is isomorphic
to one of the following Poisson algebras:
$$
\left\{
\begin{array}{ll}
e_1^2=\alpha e_3 \\
e_1 \cdot e_2=-e_2 \cdot e_1=e_2 \\
e_1 \cdot e_3=e_3 \cdot e_1=\beta e_3 \\
e_2^2=0 \\
e_2 \cdot e_3=e_3 \cdot e_2=0 \\
e_3^2=\gamma e_3
\end{array}
\right. \ \rm{with} \ \beta^2=\alpha \gamma, \ \
\left\{
\begin{array}{ll}
e_1^2=0 \\
e_1 \cdot e_2=e_2 \\
e_2 \cdot e_1=-e_2  \\
e_1 \cdot e_3=e_3 \cdot e_1=\gamma e_1 \\
e_2^2=0 \\
e_2 \cdot e_3=e_3 \cdot e_2 = \gamma e_2 \\
e_3^2=\gamma e_3
\end{array}
\right. .
$$
The first family give the Poisson algebras
$$
\mathcal{P}_{3,3}(\alpha )=
\left\{
\begin{array}{ll}
e_1^2= \alpha ^2 e_3 \\
e_1 \cdot e_2=-e_2 \cdot e_1=e_2 \\
e_1 \cdot e_3=e_3 \cdot e_1= \alpha e_3 \\
e_2^2=0 \\
e_2 \cdot e_3=e_3 \cdot e_2=0 \\
e_3^2= e_3
\end{array}
\right.  , \
\mathcal{P}_{3,4}=
\left\{
\begin{array}{ll}
e_1^2= e_3 \\
e_1 \cdot e_2=-e_2 \cdot e_1=e_2 \\
e_1 \cdot e_3=e_3 \cdot e_1=0 \qquad {\mbox {\rm and}}\\
e_2^2=e_3^2=0  \\
e_2 \cdot e_3=e_3 \cdot e_2=0
\end{array}
\right.  
$$
$$
\mathcal{P}_{3,5}=
\left\{
\begin{array}{ll}
e_1^2=e_2^2=e_3^2=0  \\
e_1 \cdot e_2=-e_2 \cdot e_1=e_2 \\
e_1 \cdot e_3=e_3 \cdot e_1=0 \\
e_2 \cdot e_3=e_3 \cdot e_2=0
\end{array}
\right. .
$$
The second family reduces to
$$
\mathcal{P}_{3,6}=
\left\{
\begin{array}{ll}
e_1^2=e_2^2=0 \\
e_3^2= e_3  \\
e_1 \cdot e_2=e_2 \\
e_2 \cdot e_1=-e_2 \\
e_1 \cdot e_3=e_3 \cdot e_1= e_1 \\
e_2 \cdot e_3=e_3 \cdot e_2= e_2
\end{array}
\right. .
$$

ii) The multiplication is given by $\{e_1,e_2\}=e_2$ and $\{e_1,e_3\}=\alpha e_3$
with $\alpha \neq 0$. From Lemma 14,
$\p$ is isomorphic to
$$
\mathcal{P}_{3,7}(\alpha)=
\left\{
\begin{array}{ll}
e_1^2=e_2^2=e_3^2=0 \\
e_1 \cdot e_2=-e_2 \cdot e_1=e_2 \\
e_1 \cdot e_3=-e_3 \cdot e_1=\alpha e_3 \\
e_2 \cdot e_3=e_3 \cdot e_2=0
\end{array}
\right.
,\ \alpha \neq 0.
$$

iii) The multiplication is given by $\{e_1,e_2\}=e_2+e_3$ and $\{e_1,e_3\}= e_3$.
As 1 is an eigenvalue of $ad_{e_1}$ with multiplicity 2,   by adapting the proof of 
Lemma~\ref{lemma2}, 
we can conclude
that $\mathcal{A}_\mathcal{P}$
is trivial. We
get the Poisson algebra:
$$
\mathcal{P}_{3,8}=
\left\{
\begin{array}{ll}
e_1^2=e_2^2=e_3^2=0 \\
e_1 \cdot e_2=-e_2 \cdot e_1=e_2+e_3 \\
e_1 \cdot e_3=-e_3 \cdot e_1= e_3 \\
e_2 \cdot e_3=e_3 \cdot e_2=0 \\
\end{array}
\right.
.
$$
$\bullet$ If $\frak{g}_\mathcal{P}$ is simple, it is isomorphic to $sl(2)$. Therefore, it is
rigid. We have already studied this case in 
the previous section. We deduce that
 $\mathcal{P}$ is isomorphic to
$$
\mathcal{P}_{3,9}=
\left\{
\begin{array}{ll}
e_1^2=e_2^2=e_3^2=0 \\
e_1 \cdot e_2=-e_2 \cdot e_1=2e_2 \\
e_1 \cdot e_3=-e_3 \cdot e_1= -2e_3 \\
e_2 \cdot e_3=-e_3 \cdot e_2=e_1 \\
\end{array}
\right.
.
$$

\section{Cohomology of Poisson algebras}

In \cite{L}, A. Lichnerowicz introduced  a cohomology for Poisson algebras. 
The $k$-cochains are skew-symmetric $k$-linear
maps that are derivatives in each of their arguments. The coboundary operator denoted by
$\delta _{LP}$ is given by

$$\delta_{LP} \varphi (X_0,...,X_{k})=
\sum_{i=0}^k (-1)^i
\left[
X_i,\varphi(X_0,...,\check{X}_i,...,X_k)\right]
+$$
$$
\sum_{0 \leq i < j \leq k }
(-1)^{i+j}\varphi(\{X_i,X_j\},X_0,...,\check{X}_i,...,\check{X}_j,...,X_k)$$
where $\check{X}_i$ means that the term $X_i$ is omitted and
$\{,\}$ is the Lie bracket of the Poisson
multiplication. Note that if $f:\mathcal{P}_1 \rightarrow
\mathcal{P}_2$ is a morphism of Poisson algebras, then $f$ does
not lead, in general, to a nontrivial functorial morphism between  the cohomology
groups. The functoriality question for 
Poisson cohomology has been addressed in the literature for instance in 
\cite{H.M}.

Since the Lichnerowicz cohomology pays attention only to the Lie part
of a Poisson algebra, we need a better definition  of cohomology that would
govern general deformations of Poisson algebras. 
Such a definition is provided by theory of quadratic Koszul operads. We describe
it in details only in degrees $0$, $1$, $2$ and $3$. Our approach 
will be based on the definition of admissible Poisson algebra.

\subsection{The operad ${\mathcal{P}}oiss$}

The operad $\mathcal{P}oiss$ has already been studied in
\cite{M.S.S}. We will give an alternative description based on the definition of Poisson algebras.
Let
$E=\Bbb{K} \left[ \Sigma_2 \right]$ be the $\Bbb{K}$-group algebra
of the symmetric group on two elements. The basis of the free
$\mathbb{K}$-module $\mathcal{F}(E)(n)$ consists of the
"parenthesized products" of $n$ variables $\left\{ x_1, ...,x_n
\right\}$. Let $R$ be the $\Bbb{K}\left[
\Sigma_{3}\right]$-submodule of $\mathcal{F(}E)(3)$ generated by
the vector
$$u=3x_1(x_2x_3)-3(x_1x_2)x_3+(x_1x_3)x_2+(x_2x_3)x_1-(x_2x_1)x_3-(x_3x_1)x_2.$$
Then $\mathcal{P}oiss$ is the binary quadratic operad with generators $E$ and relations $R$. 
It is given by
$$\mathcal{P}oiss(n)=
(\mathcal{F(}E)/\mathcal{R)}(n)=\frac{\mathcal{F(}E)(n)}{\mathcal{R(}n%
\mathcal{)}}
$$
where $\mathcal{R}$ is the operadic ideal of $\mathcal{F}(E)$ generated by
$R$ satisfying $\mathcal{R}(1)=\mathcal{R}(2)=0 $,
$\mathcal{R}(3)=R.$ The dual operad $\mathcal{P}oiss^!$ is equal
to $\mathcal{P}oiss$, that is, $\mathcal{P}oiss$ is self-dual. In
 \cite{G.R3} we defined, for a binary
quadratic operad $\mathcal{E}$, an associated quadratic operad
$\tilde{\mathcal{E}}$ which gives a functor
$$\mathcal{E} \otimes \tilde{\mathcal{E}}\rightarrow \mathcal{E}.$$ In the case $\mathcal{E}=Poiss,$ we have
$\tilde{\mathcal{E}}=\mathcal{P}oiss^!=\mathcal{P}oiss.$

\subsection{The  $k-$cochains}

We proved in \cite{G.R2} that for any
$\mathbb{K}[\Sigma_3]^2$-associative algebra
 $(\mathcal{A},\mu)$ defined by the relation
$$A_\mu^L \circ  \phi _v-A_\mu^R \circ \phi _w=0,$$
with $v,w \in \mathbb{K}[\Sigma_3]$,  the cochains $\varphi  \in
\mathcal{C}^i(\mathcal{A},\mathcal{A})$ can be chosen  invariant under
 $F_v^{\perp}\cap F_w^{\perp}$ (for the notations see
\cite{G.R2}). For a Poisson algebra we have $v=Id$, $w=3
Id-\tau_{23}+\tau_{12}-c_1+c_1^2$. Then $F_v^{\perp}\cap
F_w^{\perp}=\left\{ 0 \right\}$ and if
$\mathcal{C}^k(\mathcal{P},\mathcal{P})$ is the space of 
$k$-cochains of $\mathcal{P},$ we obtain
$$\mathcal{C}^k(\mathcal{P},\mathcal{P})=End(\mathcal{P}^{\otimes^k},\mathcal{P}). $$
{\bf Remark. } In \cite{M.S.S} an explicit presentation of the space of cochains is given using operads. 
More precisely, we have
$$\mathcal{C}^k(\mathcal{P},\mathcal{P})=
\mathcal{L}in(\mathcal{P}oiss(n)^{!}\otimes _{\Sigma _n}V^{\otimes ^n},V)$$
where $V$ is the underlying vector space (here $\C^n$). We can see that
$End(\mathcal{P}^{\otimes^k},\mathcal{P}) $ is isomorphic to
$\mathcal{L}in(\mathcal{P}oiss(n)^{!}\otimes _{\Sigma _n}V^{\otimes ^n},V)$.

\subsection{The  coboundary operators $\delta  _{\mathcal{P}}^k$, $(k=0,1,2)$}

{\bf Notation.} Let $\p$ be a Poisson algebra,
$\frak{g}_\mathcal{P}$ and $\mathcal{A}_{\mathcal{P}}$ its
corresponding Lie and associative algebras. We denote by 
$$H_C^{\star}(\frak{g}_\mathcal{P},\frak{g}_\mathcal{P})=Z_C^*(\frak{g}_\mathcal{P},\frak{g}_\mathcal{P})/
B_C^*(\frak{g}_\mathcal{P},\frak{g}_\mathcal{P})$$ the Chevalley
cohomology of $\frak{g}_{\mathcal{P}}$ and by
$H_H^{\star}(\mathcal{A}_\mathcal{P},\mathcal{A}_\mathcal{P})$ the
Harrison cohomology of $\mathcal{A}_{\mathcal{P}}$.
We will define  coboundary operators $\delta^k _{\mathcal{P}}$ on $\mathcal{C}^k(\mathcal{P},\mathcal{P})$.

\medskip

\noindent $i) \ k=0.$

\noindent We put

\medskip

\noindent  $\ H^{0}(\mathcal{P},\mathcal{P})= \left\{ X \in \mathcal{P} \, {\mbox{\rm such that}} \,
\forall Y \in \mathcal{P}, X\cdot Y=0 \right\} $.

\medskip

\noindent $ii) \ k=1.$

\noindent For $f \in End (\mathcal{P},\mathcal{P}),$ we put
$$\delta^1_\mathcal{P}f(X,Y)=f(X)\cdot Y+X \cdot f(Y)-f(X\cdot Y) $$
for any $ X,Y \in \mathcal{P}.$ Then we have

$$ H^{1}(\mathcal{P},\mathcal{P})= \,  \! H_C^{1}(\frak{g}_{\mathcal{P}},\frak{g}_{\mathcal{P}})
\cap \,
 \! H_H^{1}( \mathcal{A}_{\mathcal{P}},\mathcal{A}_{\mathcal{P}}).$$

\noindent $iii) \ k=2.$

\noindent For $\varphi \in \mathcal{C}^2(\mathcal{P},\mathcal{P})$ we define
$$
\begin{array}{ll}
\delta^2_\mathcal{P}\varphi (X,Y,Z)= & 3\varphi( X \cdot Y,Z)-3\varphi (X,Y\cdot Z)-\varphi( X \cdot Z,Y)-
\varphi( Y \cdot Z,X)
\\ & +\varphi( Y \cdot X,Z)+\varphi( Z \cdot X,Y) +3\varphi(X,Y)\cdot Z-3 X \cdot \varphi (Y,Z)
\\ & -\varphi (X,Z) \cdot Y- \varphi (Y,Z)\cdot X + \varphi (Y,X) \cdot Z
+\varphi (Z,X) \cdot Y.
\end{array}
$$
The space $ H^2(\mathcal{P},\mathcal{P})$ parametrizes 
deformations of the multiplication of $\mathcal{P}$. We saw
in the previous sections that deformations of $\p$ induce
deformations of $\frak{g}_\mathcal{P}$ and of
$\mathcal{A}_\mathcal{P}$. In constrast to $H^*(\pp,\pp)$, the 
Lichnerowicz-Poisson cohomology reflects deformations of the bracket only.

Suppose that the Poisson product satisfies $X \cdot Y=-Y\cdot X.$ Then $\{X,Y\}=X \cdot Y$ 
and $X \bullet Y=0.$
If $\varphi  \in \mathcal{C}^2(\mathcal{P},\mathcal{P})$ is also skew-symmetric, then
$$\begin{array}{rl}
\delta_\mathcal{P}^2 \varphi(X,Y,Z)= & 2\varphi(X \cdot Y,Z)+2\varphi (Y\cdot Z,X)-2\varphi (X\cdot Z,Y)\\
& +2\varphi (X,Y)\cdot Z+2\varphi (Y,Z) \cdot X-2\varphi (X,Z)\cdot Y \\
= &\delta_{LP}^2 \varphi(X,Y,Z).
\end{array}$$

\noindent We recognize the formula of Lichnerowicz-Poisson differential.

\medskip

\begin{proposition}
\label{chev}
Let $\varphi $ be in $\mathcal{C}^2(\mathcal{P},\mathcal{P})$. If $\varphi _a$ and $\varphi _s$ are respectively
the skew-symmetric and the symmetric parts of $\varphi $ then we have :

$$ \ \begin{array}{ll}
12 \delta ^2_C\varphi _a(X,Y,Z) = & \delta_\mathcal{P}^2 \varphi(X,Y,Z)-\delta_\mathcal{P}^2 \varphi(Y,X,Z)
- \delta_\mathcal{P}^2 \varphi(Z,Y,X)\\ & -\delta_\mathcal{P}^2 \varphi(X,Z,Y)
 +\delta_\mathcal{P}^2 \varphi(Y,Z,X)+\delta_\mathcal{P}^2 \varphi(Z,X,Y). \\
& \\
12\delta ^2_H\varphi _s(X,Y,Z) = & \delta_\mathcal{P}^2 \varphi(X,Y,Z)
-\delta_\mathcal{P}^2 \varphi(Z,Y,X)
+ \delta_\mathcal{P}^2 \varphi(X,Z,Y)\\
& -\delta_\mathcal{P}^2 \varphi(Z,X,Y).
\end{array}
$$
\end{proposition}
{\it Proof. } The proof  is a straightforward calculation. Recall that
 if $\varphi $ is a skew-symmetric bilinear map then the Chevalley coboundary operator 
is given by
$$
\begin{array}{lll}
\delta _C(\varphi )(X,Y,Z) &=&
\{\varphi (X ,Y),Z\}+\{\varphi (Y ,Z),X\}+ \{\varphi (Z ,X),Y\}
\\
&&
+\varphi (\{X,Y\},Z)+\varphi (\{Y,Z\},X)+\varphi (\{Z,X\},Y)\\
\end{array}
$$
and if $\varphi $ is a symmetric bilinear map then the Harrison coboundary operator 
is given by
$$
\begin{array}{lll}
\delta _H(\varphi )(X,Y,Z) &=&
\varphi (X,Y)\bullet Z-X\bullet \varphi (Y,Z)+\varphi (X\bullet Y,Z)
\\
&&
-\varphi (X,Y\bullet Z).\\
\end{array}
$$
Now, to compute
$\delta ^2_C\varphi _a$ we replace $\varphi _a(X,Y)$ by $(\varphi (X,Y)-\varphi (Y,X))/2$
and $\{X,Y\}$ by $(X \cdot Y -Y\cdot X)/2$ in the expression of $\delta ^2_C\varphi _a(X,Y,Z)$.
We leave it to the reader. 
\begin{corollary}
\label{prop}
Let $\varphi_s$ and $\varphi_a$ be the symmetric and skew-symmetric parts of
$\varphi  \in \mathcal{C}^2(\mathcal{P},\mathcal{P}).$ If
 $\varphi  \in Z^2(\mathcal{P},\mathcal{P}),$ then
$\varphi_s  \in \,  \! Z_H^2(\mathcal{A}_\mathcal{P},\mathcal{A}_\mathcal{P})$ and
$ \varphi_a  \in \,  \! Z_C^2(\frak{g}_\mathcal{P},\frak{g}_\mathcal{P}). $
\end{corollary}

\subsection{Relation between $Z^2(\mathcal{P},\mathcal{P})$ and $Z_H^2(\mathcal{A}_\mathcal{P},\mathcal{A}_\mathcal{P})$
, $Z_C^2(\frak{g}_\mathcal{P},\frak{g}_\mathcal{P})$}
To show the relation between $Z^2(\mathcal{P},\mathcal{P})$ and the classical Chevalley and Harrison
cohomological spaces, we have to introduce the following operators
$$\mathcal{L}_1, \mathcal{L}_2 : \mathcal{C}^2(\mathcal{P},\mathcal{P})\rightarrow 
\mathcal{C}^3(\mathcal{P},\mathcal{P}).$$
They are given by
$$\mathcal{L}_1(\varphi )(X,Y,Z)=\varphi (X\bullet Y,Z)-\varphi (X,Z)\bullet Y-X\bullet \varphi (Y,Z)$$
and
$$\mathcal{L}_2(\varphi )(X,Y,Z)=-3\varphi (X,\{Y,Z\})+\{\varphi (X,Y),Z\}
-\{\varphi (X,Z),Y\}.$$

\begin{lemma}
\label{lem}
Let $\varphi \in C^2(\mathcal{P},\mathcal{P})$. If $\varphi_s$ and $\varphi_a$ are the 
symmetric and skew-symmetric parts of
$\varphi$, we have
$$
\begin{array}{ll}
\delta_\mathcal{P}^2 \varphi = & \delta_C^2 \varphi_a+2\delta_H^2 \varphi_s
+\tilde \delta_C^2 \varphi_s+\tilde \delta_H^2 \varphi_a 
 + \mathcal{L}_1(\varphi_a)+ \mathcal{L}_2(\varphi_s)
\end{array}
$$
where $\tilde \delta  _C$ and $\tilde \delta _H$ are the linear maps 
$ C^2(\mathcal{P},\mathcal{P}) \rightarrow C^3(\mathcal{P},\mathcal{P})$ extending naturally
$\delta _C$ and $\delta _H$.
\end{lemma}
{\it Proof.} Starting from $\varphi =\varphi _a+\varphi _s$ and $X\cdot Y=\{X,Y\}+X\bullet Y$ we obtain
$$
\begin{array}{l}
\delta^2_\mathcal{P}\varphi (X,Y,Z)=  3\varphi_a( \{X , Y\},Z)-3\varphi_a (X,\{Y, Z\})-
\varphi_a( \{X, Z\},Y)- \\
\varphi_a( \{Y , Z\},X)
  +\varphi_a( \{Y , X\},Z)+\varphi_a( \{Z , X\},Y) +3\{\varphi_a(X,Y), Z\}\\
-3 \{X , \varphi_a (Y,Z)\}
 -\{\varphi_a (X,Z) , Y\}- \{\varphi_a (Y,Z), X\} + \{\varphi_a (Y,X) , Z\}\\
+\{\varphi_a (Z,X) , Y\}
 + 3\varphi_a( X \bullet Y,Z)-3\varphi_a (X,Y\bullet Z)-\varphi_a( X \bullet Z,Y) \\-
\varphi_a( Y \bullet Z,X)
 +\varphi_a( Y \bullet X,Z)+\varphi_a( Z \bullet X,Y) +3\varphi_a(X,Y)\bullet Z \\
-3 X \bullet \varphi_a (Y,Z)
  -\varphi_a (X,Z) \bullet Y- \varphi_a (Y,Z)\bullet X + \varphi_a (Y,X) \bullet Z \\
+\varphi_a (Z,X) \bullet Y + 3\varphi_s( \{X , Y\},Z)-3\varphi_s (X,\{Y, Z\})-
\varphi_s( \{X, Z\},Y) \\
-\varphi_s( \{Y , Z\},X)
  +\varphi_s( \{Y , X\},Z)+\varphi_s( \{Z , X\},Y) +3\{\varphi_s(X,Y), Z\}\\
-3 \{X , \varphi_s (Y,Z)\}
 -\{\varphi_s (X,Z) , Y\}- \{\varphi_s (Y,Z), X\} + \{\varphi_s (Y,X) , Z\} \\
+\{\varphi_s (Z,X) , Y\}
 + 3\varphi_s( X \bullet Y,Z)-3\varphi_s (X,Y\bullet Z)-\varphi_s( X \bullet Z,Y) \\
-\varphi_s( Y \bullet Z,X)
 +\varphi_s( Y \bullet X,Z)+\varphi_s( Z \bullet X,Y) +3\varphi_s(X,Y)\bullet Z \\-3 X \bullet \varphi_s (Y,Z)
  -\varphi_s (X,Z) \bullet Y- \varphi_s (Y,Z)\bullet X + \varphi_s (Y,X) \bullet Z \\
+\varphi_s (Z,X) \bullet Y\\
\end{array}
$$
As $\varphi _a$ is skew-symmetric and $\varphi _s$ symmetric, this relation gives
$$
\begin{array}{l}
\delta^2_\mathcal{P}\varphi (X,Y,Z)=  2\varphi_a( \{X , Y\},Z)-2\varphi_a (X,\{Y, Z\})-2
\varphi_a( \{X, Z\},Y)
\\   +2\{\varphi_a(X,Y), Z\}
-2 \{X , \varphi_a (Y,Z)\} -2\{\varphi_a (X,Z) , Y\}  
+ 4\varphi_a( X \bullet Y,Z) \\
-2\varphi_a (X,Y\bullet Z)
+2\varphi_a(X,Y)\bullet Z-4 X \bullet \varphi_a (Y,Z)
 -2\varphi_a (X,Z) \bullet Y\\
  +2\varphi_s( \{X , Y\},Z)-4\varphi_s (X,\{Y, Z\})-
2\varphi_s( \{X, Z\},Y) 
+4\{\varphi_s(X,Y), Z\} \\
-2 \{X , \varphi_s (Y,Z)\} 
 + 4\varphi_s( X \bullet Y,Z)
-4\varphi_s (X,Y\bullet Z)
  +4\varphi_s(X,Y)\bullet Z \\-4 X \bullet \varphi_s (Y,Z)\\
\end{array}
$$
that is
$$
\begin{array}{l}
\delta^2_\mathcal{P}\varphi (X,Y,Z)=  
 2\delta _C\varphi _a(X,Y,Z)+2\tilde \delta _H\varphi _a(X,Y,Z)
+2\mathcal{L}_1(\varphi_a)(X,Y,Z)\\
+4\delta _H\varphi_s(X,Y,Z)+2\tilde \delta _C\varphi _s(X,Y,Z) +2\mathcal{L}_2(\varphi_s)(X,Y,Z)
\end{array}
$$
this gives the lemma.

\begin{theorem}
\label{theo}
Let $\varphi $ be in $ C^2(\mathcal{P},\mathcal{P})$and let $\varphi_s$ , $\varphi_a$ be its 
symmetric and skew-symmetric parts. Then the following propositions are equivalent:

1. $\delta ^2_{\mathcal{P}}\varphi =0.$

\medskip

2. $\left\{
\begin{array}{l}
i) \ \delta ^2_C\varphi _a= 0, \ \delta ^2_H\varphi _s = 0. \\
ii) \ \tilde \delta_C^2 \varphi_s+\tilde \delta_H^2 \varphi_a 
 + \mathcal{L}_1(\varphi_a)+ \mathcal{L}_2(\varphi_s) =0
\end{array}
\right. $
\end{theorem} 
{\it Proof.} $2\Rightarrow 1$ is a consequence of Corollary \ref{prop}. $1\Rightarrow 2$ is a consequence of
Corollary \ref{prop} and Lemma \ref{lem}.

\medskip
\noindent{\bf Applications.} 

\noindent Suppose that $\varphi $ is skew-symmetric. Then $\varphi =\varphi _a$ and $\varphi _s=0$.
Then $\delta ^2_{\mathcal{P}}\varphi =0$ if and only if $\ \delta ^2_C\varphi = 0$
and $ \tilde \delta_H^2 \varphi 
 + \mathcal{L}_1(\varphi) =0.$ Moreless if we suppose than $\varphi $ is a biderivation 
on each argument,
 that is $\mathcal{L}_1(\varphi) =0$, then Theorem \ref{theo} implies that $\delta ^2_{\mathcal{P}}\varphi =0$
 if and only if $ \tilde \delta_H^2 \varphi =0.$ But
 $$
 \begin{array}{ll}
  \tilde \delta_H^2 \varphi (X,Y,Z)&=  \varphi (X,Y)\bullet Z-X\bullet \varphi (Y,Z)
  +\varphi (X\bullet Y,Z)-\varphi (X,Y\bullet Z)\\
  & = \mathcal{L}_1(\varphi)(X,Y,Z)+\mathcal{L}_1(\varphi)(Y,Z,X)
  \end{array}
  $$
  Thus $ \tilde \delta_H^2 \varphi =0$ as soon as $\mathcal{L}_1(\varphi) =0.$
\begin{proposition}
Let $\varphi $ be a skew-symmetric map which is a biderivation, that is $\varphi $ is a Lichnerowicz-Poisson
$2$-cochain. Then $\varphi \in Z^2_{LP}(\mathcal{P},\mathcal{P})$ if and only if 
$\varphi \in Z^2_{\mathcal{P}}(\mathcal{P},\mathcal{P})$.
\end{proposition}
Similary, if $\varphi $ is symmetric, then $\delta ^2_{\mathcal{P}}\varphi =0$ 
if and only if $\ \delta ^2_H\varphi = 0$
and $ \tilde \delta_C^2 \varphi 
 + \mathcal{L}_2(\varphi) =0.$

\subsection{The case k=3}

We need to define $\delta_\mathcal{P}^3 \psi$ for $\psi \in \mathcal{C}^3(\mathcal{P},\mathcal{P})$
so that $H^3(\mathcal{P},\mathcal{P})$ represents obstructions to 
integrability of infinitesimal deformations of the Poisson algebra
$\mathcal{P}.$ For each $\psi \in \mathcal{C}^3(\mathcal{P},\mathcal{P})$ we consider
$$
\begin{array}{ll}
\hat{\psi }(Z,T,X \cdot Y)= &\psi(Z,T,X \cdot Y)-\psi (Z,T\cdot X,Y)+\frac{1}{3}\psi (Z,T \cdot Y,X) \\
&+\frac{1}{3}\psi (Z,X \cdot Y,T)
-\frac{1}{3}\psi (Z,X \cdot T,Y)-\frac{1}{3} \psi (Z,Y \cdot T,X).
\end{array}
$$
Suppose that  $X \cdot \psi (Y,Z,T)$ appears in 
$\delta_\mathcal{P}^3 \psi(X,Y,Z,T)$. Since 
$\delta_\mathcal{P}^3 \circ \delta_\mathcal{P}^2 \varphi =0,$ we see that
the term $X\cdot \varphi (Y, Z \cdot T)$ occurs in
$X \cdot \delta_\mathcal{P}^2 \varphi (Y,Z,T).$ 
This term  appears only once
 if $\varphi $ is not skew-symmetric. Thus, in the general  case, $\delta_\mathcal{P}^3 \psi(X,Y,Z,T)$
cannot contain terms as $X \cdot  \psi (Y,Z,T).$
We conclude that $\delta_\mathcal{P}^3 \psi(X,Y,Z,T)$ can be written as:
$$\begin{array}{ll}
\delta_\mathcal{P}^3 \psi(X,Y,Z,T)=&\alpha _1 \hat{\psi}(Z,T,X \cdot Y)+\alpha _2\hat{\psi} (Y,T,X\cdot Y)+
\alpha _3\hat{\psi} (Y,Z,X \cdot T)\\
&+\alpha _4\hat{\psi} (X,T,Y \cdot Z) +\alpha _5\hat{\psi} (X,Z,Y \cdot T)+\alpha _6\hat{\psi} (X,Y,Z \cdot T).
\end{array}$$
From the relations between $Z^2(\mathcal{P},\mathcal{P})$ and 
$Z_H^2(\mathcal{A}_\mathcal{P},\mathcal{A}_\mathcal{P})$
, $Z_C^2(\frak{g}_\mathcal{P},\frak{g}_\mathcal{P})$, we have to assume that
$\delta_\mathcal{P}^3 \psi(X,Y,Z,T)=0$ as soon as $\psi $ is Lichnerowicz-¨Poisson cochain. This permits
to compute the constants $\alpha _i$. We will go in detail on this computation in a forthcomming paper.

\section{Deformations of complex Poisson algebras}

\subsection{Generalities }

By a deformation we
understand a formal deformation in Gerstenhaber's sense. It turns out that formal 
deformations are equivalent to perturbations in the sense of \cite{G.R1}. 

Let $\mathcal{P}=(V,\mu )$ be a Poisson algebra with
multiplication  $\mu$ and $V$ the underlying complex vector
space.  Let $\C[[t]]$ be the ring of complex formal power series. A
deformation of $\mu$ (or $\mathcal{P}$) is a $\C$-bilinear map:
$$\mu ' : V \times V \longrightarrow V\otimes \C[[t]]$$
given by
$$ \mu '(X,Y)=\mu(X,Y)+t\varphi_1(X,Y)+t^2\varphi_2(X,Y)+\cdots +t^n\varphi_n(X,Y)+\cdots$$
for all $X,Y \in V$ such that $\varphi_i$ are bilinear maps  satisfying, for $k \geq 1,$
$$
\left\{
\begin{array}{l}
\sum_{i+j=2k+1}\varphi_i \circ \varphi_j +\varphi_j \circ \varphi_i +\delta (\varphi_{k+1})=0,\\
\\
\sum_{i+j=2k, i<j}\varphi_i \circ \varphi_j +\varphi_j \circ \varphi_i
+\delta (\varphi_{k})+\varphi_k \circ \varphi_k=0, \\
\end{array}
\right.
$$
with
$$
\begin{array}{ll}
\varphi_i \circ \varphi_j(X,Y,Z)= & \varphi_i(\varphi_j(X,Y),Z)+)-\varphi_i(X,\varphi_j(Y,Z)) -\frac{1}{3}\varphi_i(\varphi_j(X,Z),Y)\\
& -\frac{1}{3}\varphi_i(\varphi_j(Y,Z),X)+\frac{1}{3}
\varphi_i(\varphi_j(Y,X),Z)+\frac{1}{3}\varphi_i(\varphi_j(Z,X),Y)
\end{array}
$$
and $\delta \varphi_i$ the coboundary
operator of the Poisson cohomology defined in the previous section.

\begin{definition}
A Poisson algebra $\mathcal{P}=(V,\mu )$ is rigid if every
deformation $\mu '$ is isomorphic to $\mu $, i.e., if there exists
$f \in Gl(V\otimes \C[[t]])$ such that
$$f^{ -1} (\mu (f(X),f(Y)))=\mu '(X,Y)$$
for all $X, Y \in V.$
\end{definition}
As for Lie or associative algebras, one can show, using similar arguments:
\begin{proposition}
If $H^2(\mathcal{P}, \mathcal{P})=0$, then $\mathcal{P}=(V, \mu)$  is rigid.
\end{proposition}
The converse is not true. A rigid  complex $n$-dimensional Poisson algebra with $H^2(\mathcal{P},\mathcal{P}) \neq 0$
corresponds to a point $\mu$ of the algebraic variety of Poisson structures on $\C ^n$ such that the corresponding
affine schema is not reduced at this point. We will see an example in the following section.

\subsection{Finite dimensional complex rigid Poisson algebras}

Let $\mathcal{P}=(\C ^n, \mu)$  be an $n$-dimensional complex Poisson algebra and suppose that the associated Lie algebra
$\mathfrak{g}_\mathcal{P}$ is a finite dimensional rigid solvable Lie algebra. 
It follows from \cite{A.G} that  $\mathfrak{g}_\mathcal{P}$ can be
written as $\mathfrak{g}_\mathcal{P}= \frak{t} \oplus \frak{n},$ where $\frak{n}$ is the nilradical and $\frak{t}$
a maximal abelian subalgebra such that the operators $ad X$ are semi-simple
for all $X$ in $\frak{t}$. The subalgebra
 $\frak{t}$
is called the maximal exterior torus and its
dimension  the rank of $\mathfrak{g}_\mathcal{P}$.

Suppose that $dim \, \frak{t}=1$ and for $X \in
\mathfrak{g}_\mathcal{P}$, $X \neq 0$, the restriction of the
operator $ad X$ on $\frak{n}$ is invertible (all known
solvable rigid Lie algebras satisfy this hypothesis). By Lemma 14,
the associated algebra $\mathcal{A}_\mathcal{P}$ satisfies
$\mathcal{A}_\mathcal{P}^2= \{0\}$.
\begin{theorem}
\label{theo22}
Let $\mathcal{P}$ a complex Poisson algebra such that
$\frak{g}_\mathcal{P}$ is  rigid solvable  of rank 1 (i.e $dim \, \frak{t}=1$) with
non-zero roots. Then
$\mathcal{P}$ is a rigid Poisson algebra.
\end{theorem}
{\it Proof. } If $\mu '$ is a deformation of $\mu$, then the corresponding
Lie bracket $\{ \ , \ \}_{\mu '}$ is a deformation of the Lie bracket $\{ \ , \ \}_{\mu }$ of
 $\mathfrak{g}_\mathcal{P}$. Since
$(\mathfrak{g}_\mathcal{P},\{ \ , \ \}_{\mu })$ is rigid, then $\{ \ , \ \}_{\mu ' }$  is
isomorphic to $\{ \ , \ \}_{\mu }$. If
we denote by $\mathcal{P}'=(\C ^n, \mu ')$ the deformation of $\mathcal{P}=(\C ^n, \mu)$,
then $\mathcal{A}_{\mathcal{P} '}$ satisfies also  $\mathcal{A}_{\mathcal{P} '}^2 =\{0\}$. So, $\mu '$ is isomorphic to
$\mu$ and $\mathcal{P}$ is rigid.

\medskip

\noindent Theorem \ref{theo22}  can be used to construct rigid Poisson algebras.
\begin{proposition}
Let $\frak{g}$ be a rigid solvable Lie algebra of rank 1  with non-zero
roots. Then there is only one Poisson algebra $\p$
such that $\frak{g}_\mathcal{P}=\frak{g}$.
It is defined by $$ X_i \cdot X_j =  \{X_i,X_j\} .$$
\end{proposition}

\medskip

\noindent{\bf Example.} The Poisson algebra $\mathcal{P}_{2,6}$ is rigid with
$dim H^2(\mathcal{P},\mathcal{P})=0$. In fact
$$Z^2(\mathcal{P},\mathcal{P})=\left\{ \varphi \in  \mathcal{C}^2(\mathcal{P},\mathcal{P}) , \,
\varphi(e_1,e_1)=\varphi(e_2,e_2)=0, \,
\varphi(e_1,e_2)=-\varphi(e_2,e_1) \right\}$$ and for every $f \in
End(\mathcal{P})$ we have $\delta f(e_1,e_1)=0=\delta f(e_2,e_2)$
and $\delta f(e_1,e_2)=-\delta f(e_2,e_1)=ae_1+be_2$. We observe
that $H^2_C(\mathfrak{g}_\mathcal{P} ,\mathfrak{g}_\mathcal{P} )=0
$.

\medskip

We can generalize the previous result to rigid
solvable Lie algebras $(\mathfrak{g}_\mathcal{P},\{ \ , \ \}_{\mu })$
of rank $r$. In this case the nilradical $\frak{n}$ is graded by the roots
 of $\frak{t}$ \cite{A.G}. If none of the roots is zero,
then using the  same arguments as in Lemma 14, we prove that
$\mathcal{A}_{\mathcal{P} }^2 =\{0\}$ and $\mathcal{P}$ is rigid. Then we have 

\begin{proposition}
\label{prop24}
Let $(\mathcal{P}, \mu)$  be an $n$-dimensional complex Poisson algebra such that $\mathfrak{g}_\mathcal{P}$ is  a solvable
rigid Lie algebra of rank $r$. If the roots are non-zero, then  $(\mathcal{P}, \mu)$
is rigid and $\mathcal{A}_{\mathcal{P} }^2 =\{0\}$.
\end{proposition}

\medskip

\begin{remark} \rm{\label{rig} We show how a rigid Lie algebra with 
$H^2_{C}(\mathfrak{g}_\mathcal{P},\mathfrak{g}_\mathcal{P}) \neq 0$ leads to
a rigid Poisson algebra with the same property.
Consider an admissible Poisson algebra satisfying the hypothesis of Proposition \ref{prop24}.
Thus $\mu =\{ \ , \}_{\mu}$ and if $\varphi \in Z^2(\mathcal{P}, \mathcal{P})$ is the first term of
a deformation of $\mu$, then $\varphi $ is a skew-symmetric map and $\delta \varphi (X,Y,Z)=(2/3) \delta _{C}\varphi (X,Y,Z).$
In particular, if $\mathfrak{g}_\mathcal{P}$  is rigid with
$H^2_{C}(\mathfrak{g}_\mathcal{P},\mathfrak{g}_\mathcal{P}) \neq 0$ then  $\mathcal{P}$ is rigid with
$H^2(\mathcal{P},\mathcal{P}) \neq 0$. 
This gives examples of rigid Poisson algebras with  non-trivial cohomology
based on the constructions \cite{G.A}.}
\end{remark}

\begin{remark} \rm{It may happen that a Poisson algebra $\mathcal{P}$ is rigid although 
$\mathfrak{g}_\mathcal{P}$ is not. An example is the
Poisson algebra $\mathcal{P}_{3,6}$ of Section 2.}
\end{remark}

\begin{remark} \rm{We can consider deformations of $\mathcal{P}$ which leave
the associated product of $\mathcal{A}_{\mathcal{P}}$ unchanged.
This means that $\varphi $ is a skew-bilinear map and, as in
Remark \ref{rig}, cocycles of the Poisson cohomology are also cocycles
of the Lichnerowicz-Poisson cohomology. In this case
$H^2(\mathcal{P},\mathcal{P})=H^2_{C}(\mathfrak{g}_\mathcal{P},\mathfrak{g}_\mathcal{P}).$}
\end{remark}

\subsection{The Poisson algebra $S(\frak{g})$}

Let $\frak{g}$ be a finite dimensional complex Lie algebra. We denote by $S(\frak{g})$ the
symmetric algebra on the vector space  $\frak{g}$. It is an associative commutative algebra.
Let $\left\{ e_1,...,e_n \right\}$ be a fixed basis of $\frak{g}$ and
$\{e_i,e_j\}=\sum^k_{i,j} C_{ij}^k e_k$ its structure constants. We define on $S(\frak{g})$ a
structure of Lie algebra  by
$$ P_0(p,q)=\sum_{i,j,k=1}^n C_{ij}^ke_k
\left(
\frac{ \partial p}{ \partial e_i} \frac{\partial q}{\partial  e_j}-
\frac{ \partial p}{ \partial e_j} \frac{\partial q}{\partial  e_i} \right) ,
$$
where $p=p(e_1,...,e_n)$ and $q=q(e_1,...,e_n) \in S(\frak{g})=\C[e_1,...,e_n]$. Let $p \bullet  q$ be the
ordinary associative product of the polynomials $p$ and $q$. The Lie bracket satisfies
the Leibniz rule with respect to this product. If
$$
\tilde{P_0}(p , q) =P_0(p,q)+p\bullet q
$$
then
$(S(\frak{g}),\tilde{P_0})$ is a
Poisson algebra.
This structure is usually called the linear Poisson structure on $S(\frak{g})$.

In this subsection we will be interested in deformations $\tilde{P}$  of $\tilde{P_0}$ 
on $S(\frak{g})$
which leave
the associated structure $({\cal{A}}_{S(\frak{g})},\bullet )$ unchanged.
We call such  deformations  Lie deformations of the Poisson algebra
$(S(\frak{g}), \tilde{P_0})$. Any deformation of the bracket $P_0$ can be expanded into
$$P=P_0+t\phi_1+\cdots+ t^k\phi _k+\cdots$$and the corresponding Lie deformation of $\tilde{P_0}$ is
$$\tilde{P}=\tilde{P_0}+t\phi _1+\cdots +t^k\phi _k+\cdots .$$
Then $\phi _1 \in Z^2_{L,P}((S(\frak{g}),\tilde{P_0}),(S(\frak{g}),\tilde{P_0}))$.

\medskip

Suppose now that $\frak{g}=\frak{t}\oplus \frak{n}$ is a complex solvable rigid Lie algebra.

\begin{proposition}
If $\frak{g}$ is a complex solvable rigid Lie algebra with $dim \, \frak{t} \geq 2$, then the Lie algebra
$(S(\frak{g}),P_0)$ is not rigid.
\end{proposition}
{\it Proof. } Let $\phi : S(\frak{g})\times
S(\frak{g})\longrightarrow S(\frak{g})$ be a skew-bilinear map
given by $\phi (X_1,X_2)=\alpha _{12}.1$ when $X_1,X_2 \in
\frak{t}$  and $\phi(Y_1,Y_2)=0$ when $Y_1, Y_2 \in \frak{g}$ but
$Y_1$ or $Y_2$ is not in $\frak{t}$. By the assumption,  
$\phi$ is a derivation in each argument, so 
 $\phi$ can be extended onto $S(\frak{g})$. It is easy to
see that $\phi \in Z^2_C(S(\frak{g}),S(\frak{g}))$. Since
$P_0+t\phi$ is not isomorphic to $P_0$, we have obtained a non-trivial deformation.

\begin{corollary}\cite{Pe} If $\frak{g}$ is a complex solvable rigid Lie algebra with $dim \, \frak{t} \geq 2$,
then the Poisson algebra $(S(\frak{g}),\tilde{P_0})$ is not rigid.
\end{corollary}

\medskip

Now we consider the case  $dim \, \frak{t}=1$.

\medskip
\begin{lemma}
The maximal exterior torus $\frak{t}$ is a Cartan subalgebra of $(S(\frak{g}),P_0)$.
\end{lemma}
{\it Proof. }We denote by $\{X,Y_1,...,Y_{n-1}\}$ a basis of $\frak{g}=\frak{t} \oplus \frak{n}$ adapted to this decomposition.
By definition of $\frak{t}$ we have  $\{X,Y_i\}=\lambda _iY_i$.
Then
$$
\left\{
\begin{array}{l}
P_0(X^i,X^j)=0 \ \ {\rm{for \ any \ }} i,j \\
P_0(X^i,Y_{j})=i\lambda _jX^{i-1}Y_{j} \\
P_0(X,XY_{j})=\lambda _jXY_{j} \\
P_0(X,Y_iY_j)=(\lambda _i+\lambda _j)Y_iY_j
\end{array}
\right.
$$
so that $ad_{P_0}X$ is a diagonal derivation of $S(\frak{g})$.

\medskip
We conclude that the  Lie algebra $(S(\frak{g}),P_0)$ is  graded by the eigenvalues of $ad_{P_0}X$. 
In \cite{G.A} families of rigid Lie algebras of rank $1$ were classified. 
This classification can be used to study $S(\frak{g})$
for a general rigid Lie algebra. We illustrate it on the case where
the eigenvalues of $ad_{\frak{g}}X$ are
$$1,2,...,n-1.$$
It follows from \cite{A.G} that,

- If $3 \leq n \leq 6 $ or $9 \leq n \leq 12$ then $ \frak{g}$ is not rigid.

- In the remaining cases, $\frak{g}$ is rigid.

\medskip

\noindent We consider a deformation of $\tilde{P_0}$ given as $\tilde{P}=\tilde{P_0}+t\phi_1 + ...$
with $\phi_1 \in Z^2_{L,P}((S(\frak{g},\tilde{P_0}),(S(\frak{g},\tilde{P_0}))$. It is clear that if $\phi_1(Y,Z)=0$
for every $Y,Z \in \frak{g}$ then $\phi _1 =0.$ Let $I_p$ be the Lie ideal of $S(\frak{g})$ whose elements are
polynomials of degree greater than or equal to $p$. If we denote by $S_p(\frak{g})$  the quotient Lie algebra $S(\frak{g})/I_{p+1}$, then
$S_p(\frak{g})=\C\{1\}\oplus K_p(\frak{g})$ where $K_p(\frak{g})$ 
is generated by polynomials of degree greater than or equal to
$1$. As $\tilde{\p}$ is a Lie deformation it preserves this decomposition. Thus we need to study the
Lie algebra $K_p(\frak{g})$. The Lie subalgebra generated by
 $\{X\}$ is a maximal exterior torus of $K_p(\frak{g})$. The vector $X$ is in the terminology of \cite{A.G} a
regular vector. The eigenvalues of $ad_{K_p(\frak{g})}X$ are $(1,2,...,n-1,n,...,p(n-1))$. Let $(S(X))$ be the
corresponding root system \cite{A.G}. It is easy to see that its rank is equal to $dim(\frak{n})-2$. This proves that
$K_p(\frak{g})$ is not rigid. But since we suppose that $\phi _1$ is a derivation in each argument, this implies that
$\phi_1(X,X^2)=0$ and the rank of $(S(X))$ is $dim(\frak{n})-1$. The grading of $K_p(\frak{g})$ 
by the roots
of $ad_{K_p(\frak{g})}X$ is preserved by such a deformation.

The cocycle $\phi _1$ leaves invariant each  of the eigenspaces
of adX. Let $k , \ k \leq n-1,$ be the smallest index such
 that $\phi _1$ restricted to the eigenspace associated to the eigenvalue $k$ of adX is non-zero. Then $H_k(\frak{g})$
is a non-rigid Lie algebra such that $\phi_1$ is a cocycle determined by a deformation. Conversely, let $\phi_1$ be a $2$-cocycle of the Lie
algebra $K_p(\frak{g})$ which is a derivation in each argument such that there exists $i$ with $\phi_1(Y_i,Y_{p-i}) \neq 0$.
Then we can extend $\phi _1$ to $S(\frak{g})$ to obtain a deformation of $S(\frak{g})$.

\medskip

\noindent {\bf Examples. }

\noindent 1. Let us suppose that  $\frak{g}$ is the two dimensional 
non-abelian rigid solvable Lie algebra with the bracket  defined
by $[X,Y]=Y$.
 Let $(S(\frak{g}),P_0)$ be the corresponding Poisson algebra. Then  $P_0(X,Y)=Y$.
If $P$ is a deformation of $P_0$,  since $dim(\frak{n})=1$, $P=P_0$ and $(S(\frak{g}),P_0)$ is rigid.

\medskip

\noindent 2. Let us suppose that $\frak{g}$ is the decomposable $3$-dimensional solvable Lie algebra whose brackets are
in the basis $\{X,Y_1,Y_2\}$ given by:
$$[X,Y_i]=iY_i, \ i=1,2.$$
This Lie algebra is not rigid but, as we argued in Section
2.2, there exists only one Poisson algebra structure  whose corresponding
Lie algebra is $\frak{g}$. This Poisson algebra is
${\cal{P}}_{3,7}(2)$ and it can be deformed into
${\cal{P}}_{3,7}(2+t)$. The corresponding cocycle of deformation
is given by $\phi (X,Y_2)=Y_2.$ It defines a deformation of
$(S(\frak{g}),P_0)$. The cases $n=4,5$ can be discussed in the same manner.

\medskip

\noindent 3. If $n=6$, then $\frak{g}$ is rigid. Its structure constants  are given by
$$
\left\{
\begin{array}{l}
\lbrack X, Y_i \rbrack =iY_i, \ \ i=1,...,5 \\
\lbrack Y_1, Y_i \rbrack =Y_{i+1}, \ \ i=2,3,4 \\
\lbrack Y_2, Y_3 \rbrack =Y_5.
\end{array}
\right.
$$
The Lie algebra $K_2(\frak{g})$ can be deformed using the cocycle $\phi _1(Y_1,Y_3)=Y_2^2.$ Then
$(S(\frak{g}),P_0)$ is not rigid.

\medskip

\noindent More generally, if we suppose $n > 12$, then $\frak{g}$
is rigid. Taking $\phi _1(Y_1,Y_2)=Y_1Y_2$ then the base change defined by
$Z_1=Y_1,Z_2=Y_2,Z_3=Y_3+tY_1Y_2,Z_i=[Y_1,Z_{i-1}]$ for $i \leq
n-2$ shows that the deformation of $(S(\frak{g}),P_0)$ given by
$\phi_1$ is isomorphic to a Poisson algebra which satisfies in particular
$$
\left\{
\begin{array}{l}
P(Y_1,Y_i)=Y_{i+1}, \ i=2,...,n-2 \\
P(Y_2,Y_3=Y_5+t(Y_1Y_4)+t^2(Y_1^2Y_3+Y_1Y_2^2)+t^3(Y_1^3Y_2) \\
\end{array}
\right.
$$
and that $(S(\frak{g}),P_0)$ is not rigid.

\end{document}